\documentclass[amssymb, amstex,amsmath,11pt, leqn]{article}
\usepackage[margin=1.4in]{geometry}
\usepackage{amssymb, amsmath}              
\usepackage{amsmath, amssymb, latexsym, amsfonts, url}
\usepackage{graphicx}
\usepackage{amsthm}
\usepackage{enumitem}
\usepackage{amssymb}
\usepackage{color}
\usepackage{epstopdf}

\renewcommand{\epsilon}{\varepsilon}
\theoremstyle{plain}
\newtheorem{thm}{Theorem}[section]

\newtheorem{lem}[thm]{Lemma}
\newtheorem{rem}[thm]{Remark}

\newtheorem*{theorem*}{Theorem}
\newtheorem*{proposition*}{Proposition}
\theoremstyle{definition}

\theoremstyle{remark}

\numberwithin{equation}{section}

\def\lam{\lambda}

\def\({\left(}
\def\){\right)}

\newcommand{\be}{\begin{eqnarray}}
\newcommand{\ee}{\end{eqnarray}}
\newcommand{\bea}{\begin{eqnarray*}}
\newcommand{\eea}{\end{eqnarray*}}

\newcommand{\n}{\nabla}

\def\o{\omega}

\def\a{\alpha}
\def\vp{\varphi}

\def\d{\delta}

\def\dis{\displaystyle}

\begin{document}
\begin{title}
{On the structure of almost Yamabe solitons}
\end{title}
\begin{author}
{ Seungsu Hwang \and Gabjin Yun$^*$}
\end{author}


\maketitle

\begin{abstract}
\noindent
In this paper, we study structures  of almost Yamabe solitons which are not necessarily gradient. First, we investigate conditions that both compact 
and  noncompact almost Yamabe solitons become trivial solitons which means the given vector field is a Killing vector field. Second, we show that  an almost Yamabe soliton whose  vector field  is closed  admits a local warped product structure with a one-dimensional base. This result can be considered as a generalization of a result in  \cite{c-s-z} and \cite{c-m-m}.
\vspace{.12in}


\end{abstract}

\setlength{\baselineskip}{15pt}

\section{Introduction}

A  Riemannian $n$-manifold $(M^n, g)$ together with a vector filed $X$ is called a Yamabe soliton if
there exists a constant  $\lam$  satisfying
\be 
\frac 12{\mathcal L}_Xg=(R-\lam) \, g, \label{eq1-1-1}
\ee
where  ${\mathcal L}$ is the Lie derivative and $R$ denotes the scalar curvature of $g$. Furthermore, if there exists a smooth function $f$ on $M$ such that $X = \n f$, then the triple  $(M^n, g, \n f)$   is called a gradient Yamabe soliton. 
 There has been a lot of work on the gradient Yamabe solitons for more than last decade and there is also a survey book (\cite{s-d} and references are therein). One of basic well-known facts on gradient Yamabe solitons is that the scalar curvature of any compact Yamabe soliton must be constant (\cite{d-s}, \cite{hsu}). Another interesting property on  gradient Yamabe solitons is that the metric can be basically expressed as a warped product type (\cite{c-s-z}, \cite{c-m-m}). In this paper, we consider these properties for almost Yamabe solitons. Barbosa and Ribeiro (\cite{b-r}) introduced  the notion of almost Yamabe soliton as a generalization
of Yamabe soliton. 
A  Riemannian $n$-manifold $(M^n, g)$ together with a vector filed $X$ is called an almost Yamabe soliton if
there exists a smooth function $\rho$ on $M$ which satisfies
\be 
\frac 12{\mathcal L}_Xg=(R-\rho) \, g. \label{eq1-1}
\ee
An almost Yamabe soliton $(M^n, g, X, \rho)$ is called a gradient almost Yamabe soliton if there exists a smooth function $f$ on $M$ such that $\n f = X$.  In this case, we denote the gradient almost Yamabe soliton by $(M^n, g, f, \rho)$.
 In case that $\rho$ is constant,  the triple  $(M^n, g, X)$   is just  a Yamabe soliton. 
 Many examples of almost Yamabe solitons are constructed by Barbosa  and Ribeiro (\cite{b-r}). 
 Let $X_0$ be a parallel vector field on ${\Bbb R}^{n+1}$ and let $X = X_0^\top = X_0 - \langle X_0, N\rangle N$, the tangential component of $X_0$ to ${\Bbb S}^n$, where $N$ is the outward unit normal vector field (position vector) on ${\Bbb S}^n$.
Then it is easy to see that $X$ is  the gradient of the height function $h$ and it is conformal vector field on ${\Bbb S}^n$ with conformal factor $-\langle X_0, N\rangle$ in the following sense
$$
\n_Y X = -\langle X_0, N\rangle Y
$$
for any tangent vector field $Y$ on ${\Bbb S}^n$. Thus, we can see that $({\Bbb S}^n, g_0, X, \rho = R-\frac{1}{n}\Delta h)$ is a compact gradient almost Yamabe soliton because $\Delta h = - n \langle X_0, N\rangle$.
A product manifold $({\Bbb R}\times_{\cosh t} {\Bbb S}^{n-1},  dt^2 + \cosh^2 t g_0, f = \cosh t, \rho = \sinh t + n)$ is a noncompact gradient almost Yamabe soliton (\cite{b-r}).

It is worth pointing out  (\ref{eq1-1}) means $X$ is a conformal vector field with conformal factor $R-\rho$.  The equation (\ref{eq1-1}) is equivalent to say that $X$ satisfies
$$
\langle \n_Y X, Z\rangle + \langle Y, \n_Z X\rangle = 2(R-\rho) \langle Y, Z\rangle
$$
for any vector fields $Y$ and $Z$  In particular, since $\langle {\mathcal L}_X g, g \rangle = {\rm trace}( {\mathcal L}_X g) = 2 {\rm div}X$, the function $\rho$  is given by
\be
R-\rho = \frac{1}{n}{\rm div}X,\label{eqn2024-11-29-3}
\ee
and $R-\rho$   is  called the conformal factor of $X$. Moreover, we say that a conformal vector field $X$ is {\it homothetic} (respectively, {\it Killing}) if its conformal factor $R-\rho$ is a constant function (respectively, identically zero). We say that $X$ is a nontrivial conformal vector field if it is a non-Killing vector field.
It is interesting to know whether $X$ is a Killing vector field, i.e., $R = \rho$ under some conditions.
In this direction, Barbosa and Ribeiro (\cite{b-r})  proved that for an almost Yamabe soliton $(M^n, g, X, \rho)$, $X$ is a Killing vector field if one of the following conditions holds:
(i) $M$ is compact and either $\int_M \left\{\frac12 {\rm Ric}(X, X) + (n-2)\langle \n \rho, X\rangle \right\} \le 0$ or $\int_M \langle \n \rho, X\rangle \ge 0$, or (ii) $|X| \in L^1(M)$ and either $R \ge \rho$ or $R \le \rho$ on $M$.

On the other hand, in \cite{s-m}, Seko and Maeta studies the structure of almost Yamabe solitons in Euclidean spaces. In particular, they proved that any almost Yamabe soliton $(M^n, g, X^\top, \rho)$ on a hypersurface in a Euclidean space ${\Bbb R}^{n+1}$ is contained in either a hyperplane or a sphere, where $X^\top$ is the tangential component of the position vector $X$.

One of our main results is the following.

\begin{thm}
Let $(M^n, g, X, \rho)$ be a compact almost Yamabe soliton with $n \ge 3$. 
\begin{itemize}
\item[{\rm (i)}] If either $0 \le \rho \le R$, or  $\langle X, \n\rho \rangle \le 0$, then  $X$ is a Killing vector field.
\item[{\rm (ii)}] If ${\rm Ric}(X, X) \le 0$, then $X$ is a parallel Killing vector field.
\item[{\rm (iii)}] If  $\int_M R^2 = 2 \int_M \rho^2$, then  $X$ is a Killing vector field and 
$(M, g, X)$ is an Yamabe soliton with $R=\rho=0$.
\end{itemize}
\end{thm}

In case of noncompact almost Yamabe solitons, we have the following two results.

\begin{thm} \label{thm2024-12-9-1-1}
Let $(M^n, g, X, \rho)$ be a complete non-comapct almost Yamabe soliton with $n \ge 3$. If ${\rm Ric} \le 0$ and $|X| \in L^2(M)$, then  $X$ is a parallel Killing vector field.
\end{thm}

\begin{thm} \label{lem2024-12-9-2--1-2}
Let $(M^n, g, X, \rho)$ be a complete non-comapct almost Yamabe soliton with $n \ge 3$. If $0 \le \rho \le R, \langle X, \n R\rangle  \ge 0$
and  $R+|X| \in L^2(M)$, then  $X$ is a  Killing vector field.
\end{thm}

An important particular case of a conformal vector field $X$ occurs when it satisfies the following condition
\be
 \n_Y X   = (R-\rho) Y\label{eqn2025-5-18-1}
 \ee
 for any vector field $Y$. In this case, we say $X$ is {\it closed}. A closed conformal vector field $X$ is said to be {\it parallel} if its conformal factor $R- \rho$ vanishes, and in this case, $|X|$ is constant.
 
 For an almost Yamabe soliton $(M^n, g, X, \rho)$  whose vector field   $X$ is closed, we have the following property.
 
\begin{thm}\label{thm2025-5-28-1}
Let  $(M^n, g, X, \rho)$ be an almost Yamabe soliton satisfying  (\ref{eq1-1}). If $X$ is closed, then either
\begin{itemize}
\item[{\rm (i)}] $X$ is a Killing vector field, or
\item[{\rm (ii)}]  $\dis{\left[Dd\psi + \frac{R \psi}{n-1} g - \psi {\rm Ric}\right]\left(\frac{X}{|X|}, \frac{X}{|X|}\right) = \Delta  \psi + \frac{R\psi}{n-1}}$, where $\psi = R -\rho$.
\end{itemize}
\end{thm}

Finally, we can extend one of main results in \cite{c-s-z} or \cite{c-m-m} on the metric behavior of gradient Yamabe solitons to almost Yamabe solitons with closed conformal vector fields.

\begin{thm}\label{lem2025-5-28-2}
Let  $(M^n, g, X, \rho)$ be an almost Yamabe soliton satisfying  (\ref{eq1-1}). Assume that $X$ is non-trivial(not Killing) closed vector field such that each level hypersurface $|X|^{-1}(t)$ is compact. Then $M$ is locally a  warped product $I \times_f \Sigma$, where $I$ is an interval, $\Sigma$ is a Riemannian $(n-1)$-manifold and $f = f(t)$ is a positive function on $I$ so that the metric tensor $g$ takes of the form
$$
g = dt^2 + f^2 g_\Sigma.
$$
\end{thm}

 \section{Preliminaries}
 
 
For a vector filed $W$ on $M$, let $\eta = W^\flat$, the dual $1$-form of $W$ defined by $\eta(Y) = \langle W, Y\rangle$. Then it follows from definition of Lie derivative that
$$
2\langle \n_Y W, Z\rangle  = {\mathcal L}_W g(Y, Z) + d\eta(Y, Z), \quad \forall \,Y, \, Z \in {\mathfrak X}(M).
$$
Now define a skew-symmetric $(1,1)$-tensor $\Phi:TM \to TM$ by
$$
\langle \Phi(Y), Z\rangle  = \frac{1}{2}d\eta(Y, Z),\quad   \forall \,Y, \, Z \in {\mathfrak X}(M).
$$
Thus, we can see that $X$ is a conformal vector field satisfying (\ref{eq1-1}) if and only if
\be
 \n_Y X  = (R-\rho)Y + \Phi(Y)\label{eqn2025-5-18-2}
 \ee
for any vector field $Y$. In particular, $X$ is closed if and only if $\Phi = 0$.

 \vspace{.15in}
 Let $(M^n, g, X, \rho)$ be an almost Yamabe solition satisfying (\ref{eq1-1}). Then $X$ is a conformal vector field and its conformal factor is $R - \rho$. Here we list some properties on the conformal vector field $X$ which are needed later.

\begin{thm}\label{thm2024-11-29-1}
Let $(M^n, g, X, \rho)$ be an almost Yamabe solition. Then we have the following.
\begin{itemize}
\item[$(1)$] $\n_Y X = (R-\rho)Y +  \Phi(Y)$ for any vector field $Y$.
\item[$(2)$] $ \frac12 \n |X|^2 = (R-\rho) X - \Phi(X)$.
\item[$(3)$] ${\rm Ric}(X, Y) = -(n-1)\langle Y, \n (R-\rho)\rangle  - \langle {\rm div}\Phi, Y\rangle$  for any vector field $Y$.
\item[$(4)$] $\frac12 \langle X, \nabla R\rangle = -(R-\rho)s - (n-1)\Delta (R-\rho)$.
\item[$(5)$] $\frac12 \Delta |X|^2 - |\n X|^2+ {\rm Ric}(X, X) +(n-2) \langle X, \n (R-\rho)\rangle =0$.
\item[$(6)$] ${\rm div} (\Phi(X))+ \langle {\rm div}\, \Phi , X\rangle = -|\Phi|^2$.
\item[$(7)$] $ |\n X|^2 = n(R-\rho)^2 + |\Phi|^2$.
\end{itemize}
\end{thm}
\begin{proof}
(1)  This is just the equation (\ref{eqn2025-5-18-2}). 

\noindent
(2) It is  easy to see that this  comes from (1).

\noindent
(3) One can find this identity in  \cite{r-u} or \cite{h-y}.

\noindent
(4) See \cite{y}.

\noindent
(5) It is due to Petersen and Wylie \cite{p-w} that, for any vector field $X$ (not necessarily conformal),
$$
{\rm div}({\mathcal L}_X g) (X) = \frac12 \Delta |X|^2 - |\n X|^2+ {\rm Ric}(X, X) + \n_X {\rm div}X.
$$
From (\ref{eq1-1}), since ${\rm div}({\mathcal L}_X g) (X)= 2\langle X, \n (R-\rho)\rangle$ and  ${\rm div}(X) = n(R- \rho)$, we  obtain 
\be
\frac12 \Delta |X|^2 - |\n X|^2+ {\rm Ric}(X, X) +(n-2) \langle X, \n (R-\rho)\rangle =0.\label{eqn2024-11-29-2}
\ee
(6) This is  proved in \cite{h-y}.

\noindent
(7) Taking the divergence of both sides in (2) above, we obtain
\be
\frac12 \Delta |X|^2  =  \langle X, \n (R-\rho)\rangle + n(R-\rho)^2 - {\rm div}(\Phi(X)).\label{eqn2023-3-14-7-1}
\ee
From (3) above and (\ref{eqn2024-11-29-2}), we have
$$
\frac12 \Delta |X|^2 = |\n X|^2 + \langle X, \n (R-\rho)\rangle + \langle {\rm div}\Phi, X\rangle.
$$
So,
$$
 |\n X|^2 = n(R-\rho)^2 - {\rm div}(\Phi(X)) - \langle {\rm div}\Phi, X\rangle.
 $$
 Finally, by (6), we obtain
 $$
  |\n X|^2 = n(R-\rho)^2 + |\Phi|^2.
  $$
\end{proof}

\section{Compact Almost Yamabe Solitons}

Assume that $(M^n, g)$ is a compact Riemannian manifold (without boundary). Applying the divergence theorem to (\ref{eqn2024-11-29-3}), we obtain
\be
\int_M (R- \rho) dV = 0.\label{eqn2024-11-29-4}
\ee
Since ${\rm div}(RX) = \langle \n R, X\rangle + R {\rm div}(X)$, by  (\ref{eqn2024-11-29-3}) again, we have
$$
\int_M R(R-\rho) = \frac{1}{n}\int_M R \, {\rm div}X = - \frac{1}{n}\int_M \langle \n R, X\rangle.
$$
On the other hand, by Theorem~\ref{thm2024-11-29-1} (4), since $\dis{\int_M R(R-\rho) = - \frac{1}{2}\int_M \langle X, \n R \rangle}$,  we have
\be
\int_M R(R-\rho)dV = 0\label{eqn2024-11-29-5}
\ee
when $n \ge 3$.

Thus, we obtain the following.

\begin{thm}
Let $(M^n, g, X, \rho)$ be a compact almost Yamabe soliton with $n \ge 3$. If $0 \le \rho \le R$, then $X$ is a Killing vector field.
\end{thm}
\begin{proof}
From (\ref{eqn2024-11-29-5}) together with our assumption, we have
$$
\int_M (R- \rho)^2 dV = - \int_M \rho (R-\rho)dV \le 0.
$$
\end{proof}

\begin{thm}
Let $(M^n, g, X, \rho)$ be a compact almost Yamabe soliton with $n \ge 3$. If $\langle X, \n\rho \rangle \le 0$, then $X$ is a Killing vector field.
\end{thm}
\begin{proof}
By Theorem~\ref{thm2024-11-29-1} (4) together with (\ref{eqn2024-11-29-5}), we have
$$
\int_M \langle X, \n R\rangle =0.
$$
So, it follows from (\ref{eqn2023-3-14-7-1}) that
$$
n\int_M (R-\rho)^2 dv_g = \int_M \langle X, \n\rho\rangle - \int_M \langle X, \n R\rangle = \int_M \langle X, \n\rho\rangle.
$$
Thus, our assumption  implies that $R= \rho$ and so $X$ is  a Killing vector field.
In this case, we have $\langle X, \n \rho\rangle = \langle X, \n R\rangle = 0$.
\end{proof}

\begin{thm}\label{lem2025-1-23-1}
Let $(M^n, g, X, \rho)$ be a compact almost Yamabe soliton with $n \ge 3$. If ${\rm Ric}(X, X) \le 0$, then  $X$ is a  Killing vector field and parallel.
\end{thm}
\begin{proof}
Since ${\rm div}[(R-\rho)X] = \langle X, \n (R-\rho)\rangle + n (R- \rho)^2$,  it follows from  Theorem~\ref{thm2024-11-29-1} (5) that
$$
\int_M \left[|\n X|^2 -{\rm Ric}(X, X)  + n(n-2)(R-\rho)^2 \right] dV = 0.
$$
This together with our assumption shows  that $R= \rho$ and so $X$ is  a parallel Killing vector field. 
\end{proof}

\begin{rem}
In Theorem~\ref{lem2025-1-23-1}, we can see that $X$ is closed in the sense of $\Phi = 0$ by  Theorem~\ref{thm2024-11-29-1} (7).
\end{rem}

\begin{thm}
Let $(M^n, g, X, \rho)$ be a compact almost Yamabe soliton with $n \ge 3$. If $\int_M R^2 = 2 \int_M \rho^2$, then
$X$ is a Killing vector field and $(M, g, X)$ is an Yamabe soliton.
\end{thm}
\begin{proof}
Let $\psi := R- \rho$. Recall that 
$$
\int_M(R-\rho) = 0 = \int_M R(R-\rho) = \int_M R \psi =0
$$
and so,
$$
\int_M (R-\rho)^2  = -\int_M \rho \psi \ge 0. 
$$
In particular, we have
\be
\int_M R\rho \le \int_M \rho^2.\label{eqn2025-1-25-5}
\ee 
Now,  by Cauchy-Schwarz inequality,  
\be
\int_M \psi^2 = \int_M (R-\rho)^2  = -\int_M \rho \psi \le \left(\int_M \rho^2\right)^{\frac{1}{2}}\left(\int_M \psi^2\right)^{\frac{1}{2}},\label{eqn2025-1-24-6}
\ee
which shows that
$$
\int_M (R^2 - 2R\rho + \rho^2) = \int_M \psi^2 \le \int_M \rho^2.
$$
Hence
$$
0 \le \int_M R^2 \le 2 \int_M R\rho \le 2 \int_M \rho^2.
$$
However, since $\int_M R^2 = 2 \int_M \rho^2$ from our assumption,  we obtain
\be
\int_M \psi^2 = \int_M \rho^2 =\frac12 \int_M R^2.\label{eqn2025-1-24-7}
\ee
This together with  (\ref{eqn2025-1-24-6}) shows  that
two functions $\psi$ and $\rho$ are parallel. Consequently, we have
$\psi = R-\rho = \rho$  or $\psi = - \rho$. Therefore, $\rho = R = 0$ by (\ref{eqn2025-1-25-5}) for the first case, and it is trivial in the second case.

\end{proof}

\section{Non-compact Almost Yamabe Solitons}

In this section, we consider a complete noncompact almost Yamabe solitons with $n \ge 3$ satisfying
(\ref{eq1-1}).

\begin{thm} \label{lem2024-12-9-1}
Let $(M^n, g, X, \rho)$ be a complete non-comapct almost Yamabe soliton with $n \ge 3$. If ${\rm Ric} \le 0$ and $|X| \in L^2(M)$, then  $X$ is a parallel Killing vector field.
\end{thm}
\begin{proof}
Choose a cut-off function $\xi$ satisfying
$$
0 \le \xi \le 1, \quad {\rm supp}(\xi) \subset B(r) \quad\mbox{and}\quad  |\n \xi|^2 \le \frac{2}{r^2},
$$
where $B(r)$ is a geodesic $r$-ball in $(M, g)$. Note that, from (\ref{eqn2024-11-29-3}),
$$
{\rm div}[\xi^2(R-\rho)X] =
\xi^2 \langle X, \n(R-\rho)\rangle + 2\xi (R-\rho)\langle X, \n \xi\rangle + n\xi^2 (R-\rho)^2
$$
and so
\be
\int_M \xi^2 \langle X, \n(R-\rho)\rangle = - 2\int_M \xi (R-\rho)\langle X, \n \xi\rangle 
- n \int_M \xi^2 (R-\rho)^2.\label{eqn2024-12-8-1}
\ee
Recall  Theorem~\ref{thm2024-11-29-1} (5):
$$
\frac12 \Delta |X|^2 - |\n X|^2+ {\rm Ric}(X, X) +(n-2) \langle X, \n (R-\rho)\rangle =0.
$$
Multiplying this  by $\xi^2$ and integrating it over $M$, we have
\bea
\int_M \xi^2 |\n X|^2 
&=&
(n-2)\int_M \xi^2 \langle X, \n (R-\rho)\rangle + \int_M \xi^2 {\rm Ric}(X, X) + \frac12 \int_M \xi^2 \Delta |X|^2\\
&=&
-n(n-2) \int_M \xi^2 (R-\rho)^2  - 2(n-2) \int_M \xi (R-\rho) \langle \n \xi, X\rangle\\
&&\qquad
+ \int_M \xi^2 {\rm Ric}(X, X) - \int_M \xi \langle \n \xi, \n |X|^2 \rangle \\
&\le& 
-n(n-2) \int_M \xi^2 (R-\rho)^2  - 2(n-2) \int_M \xi (R-\rho) \langle \n \xi, X\rangle\\
&&\qquad
+\frac12\int_M \xi^2 |\n |X||^2 +2\int_M |\n \xi|^2 |X|^2\\
&\le&
-n(n-2) \int_M \xi^2 (R-\rho)^2  +(n-2) \int_M \xi^2 (R-\rho)^2 \\
&&\qquad
+ (n-2)\int_M |\n \xi|^2 |X|^2 +\frac12\int_M \xi^2 |\n |X||^2 +2\int_M |\n \xi|^2 |X|^2.
\eea
In other words, we obtain
\bea
\int_M \xi^2 |\n X|^2 +(n-1)(n-2)\int_M \xi^2 (R-\rho)^2 \le 2n \int_M |\n \xi|^2 |X|^2 \le \frac{4n}{r^2}\int_M |X|^2.
\eea
By letting $r\to \infty$ and using the $L^2$-finiteness of $|X|$, we must have $\n X = 0$.

Finally, from (7) in Theorem~\ref{thm2024-11-29-1}, we have $R = \rho$ and $\Phi = 0$, 
which means  $X$ is a closed Killing vector field.

\end{proof}

Similarly as Theorem~\ref{lem2024-12-9-1}, we have the follwoing.

\begin{thm} \label{lem2024-12-9-2}
Let $(M^n, g, X, \rho)$ be a complete non-comapct almost Yamabe soliton with $n \ge 3$. If $0 \le \rho \le R, \langle X, \n R\rangle  \ge 0$
and  $R+|X| \in L^2(M)$, then  $X$ is a  Killing vector field.
\end{thm}
\begin{proof}
As in the proof of Theorem~\ref{lem2024-12-9-1}, choose a cut-off function $\xi$ satisfying
$$
0 \le \xi \le 1, \quad {\rm supp}(\xi) \subset B(r) \quad\mbox{and}\quad  |\n \xi|^2 \le \frac{2}{r^2}.
$$
From (\ref{eqn2024-12-8-1}), we  have
\be
n \int_M \xi^2 (R-\rho)^2  +\int_M \xi^2 \langle X, \n R\rangle =\int_M \xi^2 \langle X, \n \rho\rangle
- 2\int_M \xi (R-\rho)\langle X, \n \xi\rangle. \label{eqn2024-12-9-2}
\ee
Since
\bea
 \xi^2 \rho (R-\rho) &=& \frac{1}{n} \xi^2 \rho {\rm div}X \\
&=&
\frac{1}{n}  {\rm div} (\xi^2 \rho X)- \frac{2}{n} \xi \rho \langle X, \n \xi\rangle - \frac{1}{n} \xi^2 \langle X, \n \rho\rangle,
 \eea
 we have
 \be
 \int_M \xi^2 \langle X, \n \rho\rangle + n \int_M \xi^2 \rho (R-\rho) = -2 \int_M \xi \rho \langle X, \n \xi\rangle.\label{eqn2025-12-10-1}
 \ee
 Since $R \in L^2$ and $0 \le \rho \le R$, the function $\rho$ is also in $L^2(M)$, and so
 \bea
 \int_M \xi \rho \langle X, \n \xi\rangle &\le& \left(\int_M \xi^2 \rho^2 \right)^{\frac12} \left(\int_M |X|^2 |\n \xi|^2 \right)^{\frac12}\\
 &\le&\frac{C}{r}\left(\int_M |X|^2\right)^{\frac12}.
 \eea
 Thus, by letting $r \to \infty$, we have $\left|  \int_M  \xi \rho \langle X, \n \xi\rangle\right|  \to 0$,
  which means $\dis{\lim_{r\to \infty} \int_M \xi^2 \langle X, \n \rho \rangle =0}$ because of our assumption $0 \le \rho \le R$ 
  and (\ref{eqn2025-12-10-1}).
  
  \vspace{.12in}
 Next, from Young's inequality, we have
 $$
  2\int_M \xi (R-\rho)\langle X, \n \xi\rangle \le  \int_M \xi^2 (R-\rho)^2 + \int_M |X|^2 |\n \xi|^2.
  $$
  So, from (\ref{eqn2024-12-9-2}), we have
  \bea
  (n-1)\int_M \xi^2 (R-\rho)^2 +\int_M \xi^2 \langle X, \n R\rangle &\le&\int_M \xi^2 \langle X, \n \rho\rangle +
  \int_M |X|^2 |\n \xi|^2\\
  &\le&
  \int_M \xi^2 \langle X, \n \rho\rangle + \frac{2}{r^2}\int_M |X|^2.
 \eea
 Therefore, by letting $r\to \infty$ we must have $R - \rho =0$ from our assumptions.
 
\end{proof}

\section{Structures of almost Yamabe  solitons}

Let  $(M^n, g, X, \rho)$ be an almost Yamabe soliton satisfying (\ref{eq1-1}).
Through out this section, we may assume $X$ is {\it closed}  which implies $$\n_YX = \psi Y$$ for any vector field    $Y$, where $\psi=R - \rho$. In particular, since $X$ is closed, we have
$$
\langle \n_Y X, Z\rangle  = \langle Y, \n_Z X\rangle
$$
for any vector fields $Y$ and $Z$, and so defining $\n X(Y, Z) = \langle \n_Y X, Z \rangle$, we have $\n X = \psi g$. Moreover, since
$X$ is closed, we have $\Phi = 0$ in Theorem~\ref{thm2024-11-29-1}. 

\vspace{.2in}
\noindent
Defining $u=|X|^2$ and recalling $\psi = R -\rho$, we have
$$
\n  u = 2\psi X\quad \mbox{and}\quad |\n u|^2 = 4 \psi^2 u
$$
because $\n_Y X = \psi Y$ for any vector field $Y$. Then
\bea
\frac12 \Delta |\n u|^2 &=&2 (\Delta \psi^2 ) u + 2\psi^2 \Delta u + 4 \langle \n \psi^2, \n u\rangle\\
&=&
4\left(\psi \Delta \psi + |\n \psi|^2 \right)u + 2\psi^2 \Delta u + 8 \psi \langle \n \psi, \n u\rangle.
\eea
Recall, by Theorem~\ref{thm2024-11-29-1}, that
\bea
\Delta \psi = - \frac{R}{n-1}\psi - \frac{1}{2(n-1)}\langle X, \n R \rangle.
\eea
Also note that $\Delta u  = 2\langle \n \psi, X\rangle + 2n \psi^2$ and $\langle \n \psi, \n u\rangle = 2\psi \langle \n \psi, X\rangle$.
Thus,
\be
\frac12 \Delta |\n u|^2 = 4\left(\psi \Delta \psi + |\n \psi|^2 \right) u +20 \psi^2 \langle \n \psi, X\rangle + 4n \psi^4.\label{eqn2025-1-19-1}
\ee

\begin{lem}\label{lem2025-1-19-3}
Let  $(M^n, g, X, \rho)$ be an almost Yamabe soliton satisfying  (\ref{eq1-1}). Then
$$
\frac12 \Delta |\n u|^2 = 4 |\n \psi|^2  u + 4n \psi^4 + 4\psi Dd\psi(X, X)  +4(n+4)\psi^2  \langle \n \psi, X\rangle .
$$
\end{lem}
\begin{proof}
From Bochner formula, 
\be
\frac12 \Delta |\n u|^2 =|Ddu|^2 + \langle \n \Delta u, \n u\rangle + {\rm Ric}(\n u, \n u).\label{eqn2025-1-19-2}
\ee
Since $Ddu = 2 d\psi \otimes X + 2\psi \n X = 2 d\psi \otimes X +2n \psi^2 g$, we have
$$
|Ddu|^2 = 4 |\n \psi|^2 u + 4n \psi^2 + 8\psi^2 \langle \n \psi, X\rangle.
$$
Next, from $\Delta u  = 2\langle \n \psi, X\rangle + 2n \psi^2$ and $\n X = \psi g$, 
\bea
\n \Delta u &=& 2 Dd\psi (X) + 2  \n_{\n \psi} X + 4n \psi \n \psi\\
&=&
 2 Dd\psi (X) + 2 (2n+1) \psi \n \psi
 \eea
 and so
 $$
\langle \n \Delta u, \n u\rangle =\langle \n \Delta u, 2\psi X\rangle = 4\psi Dd\psi(X, X) + 4(2n+1) \psi^2 \langle \n \psi, X\rangle.
$$
Finally, from Theorem~\ref{thm2024-11-29-1},
$$
{\rm Ric}(\n u, \n u) = 4\psi^2 {\rm Ric}(X, X) = - 4 (n-1)\psi^2 \langle \n \psi, X\rangle.
$$
Substituting these into (\ref{eqn2025-1-19-2}), we obtain
$$
\frac12 \Delta |\n u|^2 = 4|\n \psi|^2 u + 4n \psi^4 + 4 \psi Dd\psi(X, X) + 4(n+4) \psi^2 \langle \n \psi, X\rangle.
$$

\end{proof}

\begin{lem}\label{lem2025-1-19-5}
Let  $(M^n, g, X)$ be an almost Yamabe soliton satisfying  (\ref{eq1-1}). Then
$$
\psi \left[Dd\psi + \frac{R \psi}{n-1} g - \psi {\rm Ric}\right](X, X) = \psi \left(\Delta  \psi + \frac{R\psi}{n-1}\right) |X|^2.
$$
\end{lem}
\begin{proof}
By Lemma~\ref{lem2025-1-19-3} and the equation (\ref{eqn2025-1-19-1}), we have
$$
\psi Dd\psi (X, X) - u \psi  \Delta \psi + (n-1) \psi^2 \langle \n \psi, X \rangle =0.
$$
Since ${\rm Ric}(X, X) = - (n-1)\langle \n \psi, X\rangle$ and $u = |X|^2 = g(X, X)$, we obtain
$$
\psi \left[Dd\psi + \frac{R \psi}{n-1} g - \psi {\rm Ric}\right](X, X) = \psi \left(\Delta  \psi + \frac{R\psi}{n-1}\right) |X|^2.
$$

\end{proof}

From Lemma~\ref{lem2025-1-19-5}, we have the following.

\begin{thm}\label{thm2025-1-24-1}
Let  $(M^n, g, X, \rho)$ be an almost Yamabe soliton satisfying  (\ref{eq1-1}). If $X$ is closed, then either
\begin{itemize}
\item[{\rm (i)}] $X$ is a Killing vector field, or
\item[\rm{(ii)}]  $\dis{\left[Dd\psi + \frac{R \psi}{n-1} g - \psi {\rm Ric}\right]\left(\frac{X}{|X|}, \frac{X}{|X|}\right) = \Delta  \psi + \frac{R\psi}{n-1}.}
$
\end{itemize}
\end{thm}

\vspace{.12in}
\begin{rem}
{\rm We can also prove Theorem~\ref{thm2025-1-24-1} as follows. Define $u = |X|^2$ as above.
Since $\n u = 2 \psi X$   and $\n X = \psi g$, we have
$$
Ddu = 2d\psi \otimes X^\flat + 2\psi^2 g
$$
and so
\bea
\n |\n u|^2 &=& 2 Ddu(\n u) = 4\langle\n \psi, \n u\rangle X + 4 \psi^2 \n u\\
&=&
8 \psi \langle \n \psi, X\rangle X + 4 \psi^2 \n u.
\eea
On the other hand, since $|\n u|^2 = 4 \psi^2 u$, we have
$$
\n |\n u|^2 = 8 u \psi \n \psi + 4\psi^2 \n u.
$$
Comparing these two equations, we obtain
\be
\psi u \n \psi = \psi \langle \n \psi, X\rangle X.\label{eqn2025-1-25-2}
\ee
Thus, if $X$ is not a Killing vector field, we have 
$$
u \n \psi = \langle \n \psi, X\rangle X.
$$
Taking the covariant derivative, we have
$$
du \otimes d \psi + u Dd \psi = Dd\psi(X) \otimes X^\flat + \psi d\psi\otimes X^\flat + \psi \langle \n \psi, X \rangle g.
$$
Taking the trace of both sides, we obtain
$$
\langle \n u, \n \psi\rangle + u \Delta \psi = Dd\psi(X, X) + (n+1) \psi \langle \n \psi, X \rangle.
$$
This shows immediately
\be
|X|^2 \Delta \psi = Dd\psi(X, X) -  \psi {\rm Ric}(X, X)\label{eqn2025-1-24-3-1}
\ee
becuase $\langle \n u, \n \psi\rangle = 2\psi \langle \n \psi, X\rangle, u = |X|^2$ and ${\rm Ric}(X, X) = -(n-1)\langle \n \psi, X\rangle$.
}

 \end{rem}

Let  $(M^n, g, X, \rho)$ be an  almost Yamabe soliton satisfying  (\ref{eq1-1}) and assume that $|X|^{-1}(c)$ be a hypersurface.
Then $X$ is a normal vector field to $|X|^{-1}(c)$. In fact, let $\a(t)$ be a curve in $|X|^{-1}(c)$ so that $|X|^2(\a(t)) = c^2$.
Taking the derivative with respect to $t$ and using Theorem~\ref{thm2024-11-29-1} together with $\Phi = 0$, we have
$$
2\langle \n_{\a'(t)}X, X\rangle = 2\psi \langle \a'(t), X\rangle = 0,
$$
which shows that $X$ is normal to $|X|^{-1}(c)$.

In case of compact almost Yamabe solitons $(M^n, g, X, \rho)$ or the norm of closed conformal vector field $X$ has compact level hypersurfaces, we have the following property.

\begin{lem}\label{lem2025-1-19-8}
Let  $(M^n, g, X, \rho)$ be an  almost Yamabe soliton satisfying  (\ref{eq1-1}). Assume that $X$ is non-trivial closed vector field such that each level hypersurface $|X|^{-1}(t)$ is compact. Then the function $\psi =R-\rho$ is constant along  the hypersurface $|X| = c${\rm ($\ne 0$, constant)}.
\end{lem}
\begin{proof}
Let $\Sigma$ be the hypersurface given by $ |X|=c$ and let $\{e_1, e_2, \cdots, e_{n-1}, e_n\}$ be a local frame around a point on $\Sigma$ with $e_n = \frac{X}{|X|}$ which is unit normal vector field on $\Sigma$. From Lemma~\ref{lem2025-1-19-5} together with our assumptions, we have
$$
\left[Dd\psi + \frac{R \psi}{n-1} g - \psi {\rm Ric}\right](e_n, e_n) = \Delta  \psi + \frac{R\psi}{n-1}
$$
and hence
\be
\sum_{i=1}^{n-1}\left(Dd\psi + \frac{R \psi}{n-1} g - \psi {\rm Ric}\right)(e_i, e_i) =0.\label{eqn2025-1-19-6}
\ee
Now, note that
\bea
\sum_{i=1}^{n-1}Dd\psi (e_i, e_i) &=& e_ie_i(\psi) - \n_{e_i}e_i (\psi)\\
&=&
 e_ie_i(\psi) - \left(\n_{e_i}e_i\right)^\top (\psi) -\langle \n_{e_i}e_i, e_n\rangle e_n(\psi)\\
 &=&
 \Delta^\Sigma \psi + (n-1)\frac{\psi}{|X|^2} \langle X, \n \psi\rangle.
 \eea
 Substituting this into (\ref{eqn2025-1-19-6}), we obtain
 \bea
  \Delta^\Sigma \psi + (n-1)\frac{\psi}{|X|^2} \langle X, \n \psi\rangle 
  &=& 
  -\sum_{i=1}^{n-1}\left(\frac{R \psi}{n-1} g - \psi {\rm Ric}\right)(e_i, e_i)\\
  &=&
 -\psi{\rm Ric}(e_n, e_n) = - \frac{\psi}{|X|^2} {\rm Ric}(X, X)\\
 &=&
 (n-1)\frac{\psi}{|X|^2}\langle \n \psi, X\rangle.
 \eea
 Therefore, $\Delta^\Sigma \psi = 0$, which means $\psi$ is constant on $\Sigma$ when $\Sigma$ is compact.
\end{proof}

\begin{lem}\label{lem2025-1-23-2}
Let  $(M^n, g, X, \rho)$ be an almost Yamabe soliton satisfying  (\ref{eq1-1}). Assume that $X$ is non-trivial closed vector field such that each level hypersurface $|X|^{-1}(t)$ is compact. Then each level hypersurface $|X| = c${ \rm ($\ne 0$, constant)} is totally umbilic and so has a constant mean curvature in $M$.
\end{lem}
\begin{proof}
Let $\Sigma$ be the hypersurface given by $ |X|=c$ and let $\{e_1, e_2, \cdots, e_{n-1}, e_n\}$ be a local frame around a point on $\Sigma$ with $e_n = \frac{X}{|X|}$. Since $|X|$ is constant on $\Sigma$, for each $i = 1, 2, \cdots, n-1$, we have
$$
\n_{e_i} e_n =\n_{e_i} \frac{X}{|X|} = \frac{1}{|X|} \n_{e_i}X = \frac{\psi}{|X|}e_i.
$$
By Lemma~\ref{lem2025-1-19-8},  this shows that the hypersurface $\Sigma: |X|=c$ is totally umbilic.
In particular, denoting by $\mathfrak m$ the (unnormalized) mean curvature of $\Sigma$, we have
$$
\mathfrak m = -\sum_{i=1}^{n-1}\langle \n_{e_i}e_n, e_i\rangle = -\frac{(n-1)\psi}{|X|}
$$
and so $\mathfrak m$ is constant on $\Sigma$.
\end{proof}

\begin{thm}\label{thm2025-1-25-3}
Let  $(M^n, g, X, \rho)$ be a compact almost Yamabe soliton satisfying  (\ref{eq1-1}) and assume $X$ is closed.
If either $X$ is nowhere vanishing or  $\langle X, \n \psi\rangle =0$ with $\psi = R-\rho$, then $R=\rho$ and $X$ is a Killing vector field.
\end{thm}
\begin{proof}
 Let $\{e_1, e_2, \cdots, e_{n-1}, e_n\}$ be a local frame  with $e_n = \frac{X}{|X|}$. 
Since $\n_{e_n}X = \psi e_n$, we have 
$$
e_n(|X|) = \psi\quad \mbox{and}\quad e_n \left(\frac{1}{|X|}\right) = - \frac{\psi}{|X|^2}.
$$
In particular,
\bea
\n_{e_n}e_n &=& \n_{e_n}\left(\frac{X}{|X|}\right) = e_n \left(\frac{1}{|X|}\right)  X + \frac{1}{|X|}\n_{e_n}X\\
&=&
- \frac{\psi}{|X|^2}X + \frac{\psi}{|X|} e_n\\
&=& 0.
\eea
Thus, from our assumption,
$$
Dd\psi(e_n, e_n) = e_n \langle e_n, \n \psi\rangle =0.
$$
Since ${\rm Ric}^M(X, X) = -(n-1)\langle  X, \n \psi,\rangle$, from (\ref{eqn2025-1-24-3-1}), 
$$
\Delta \psi = Dd\psi(e_n, e_n) -  \psi {\rm Ric}^M_{nn}  = (n-1) \psi \frac{\langle  X, \n \psi\rangle}{|X|^2}.
$$
Defining $u = |X|^2$, this can be written in the following form:
\be
u \Delta \psi =  (n-1) \psi \langle  X, \n \psi\rangle.\label{eqn2025-1-25-1}
\ee
Since $\n u = 2\psi X$, this equation is equivalent to
\be
\Delta \psi - \frac{n-1}{2} \langle \n \ln u, \n \psi\rangle =0.\label{eqn2025-9-11-1}
\ee
If $X$ is nowhere vanishing or $\langle X, \n \psi\rangle =0$, we can apply  the maximum principle to (\ref{eqn2025-9-11-1}) and can conclude that  the function  $\psi$ must be constant. Finally, since $\int_M \psi =\frac{1}{n}\int_M {\rm div}X = 0$, $\psi $ must be vanishing.
\end{proof}

\begin{rem}
{\rm In Theorem~\ref{thm2025-1-25-3}, the condition that $X$ is nowhere vanishing is crucial. In fact, one can construct a closed (gradient)
non-Killing vector field $X = \n h$  on ${\Bbb S}^n$ so that $({\Bbb S}^n, g_0, X, \rho)$ is an almost Yamabe soliton with $\rho = R - \frac{1}{n}\Delta h$, where $R$ is the scalar curvature  of $({\Bbb S}^n, g_0)$. See \cite{b-r} for more detail.

On the other hand, for any  any compact gradient Yamabe soliton $(M^n, g, X= \n f)$, the condition $\langle X, \n \psi\rangle = 0$ is obvious because  it  has constant scalar curvature (\cite{d-s}, \cite{hsu}). So, this condition is meaningful for almost Yamabe solitons since $\rho$ and $\psi = R-\rho$ is not constant in general.
}
\end{rem}

A smooth vector field $V$ on a Riemannian manifold $(M, g)$ is said to be {\it torqued vector field}   \cite{chen} if it satisfies 
 \be
 \n_YV = \vp Y + \o(Y)V\label{eqn2024-5-29-1}
 \ee
 for any vector field $Y$ satifying $\o(V) = 0$. Here the function $\vp$ and the $1$-form $\o$ are  called the torqued function  and  the torqued form of $V$, respectively.
 In \cite{chen} and \cite{chen2}, B-Y. Chen studied  torqued vector fields systematically. One of main results is  that  if a   Riemannian manifold $(M, g)$ admits a torqued vector field $V$, then $M$ is locally a twisted product $I \times_f \Sigma$ such that $V$ is always tangent to $I$, where $I$ is an open interval.

\begin{thm}\label{lem2025-5-18-7}
Let  $(M^n, g, X, \rho)$ be an almost Yamabe soliton satisfying  (\ref{eq1-1-1}). Assume that $X$ is non-trivial  closed vector field such that each level hypersurface $|X|^{-1}(t)$ is compact. Then $M$ is locally a  warped product $I \times_f \Sigma$, where $I$ is an interval, $\Sigma$ is a compact Riemannian $(n-1)$-manifold and $f = f(t)$ is a positive function on $I$ so that the metric tensor $g$ takes of the form
$$
g = dt^2 + f^2 g_\Sigma.
$$
In particular, if $M$ is locally product $M = I \times \Sigma$ with constant $f$, then $X$ is a Killing vector field.
\end{thm}
\begin{proof}
Let $\Sigma$ be the (compact) hypersurface given by $ |X|=c$ and let $\{e_1, e_2, \cdots, e_{n-1}, e_n\}$ be a local frame around a point on $\Sigma$ with $e_n = \frac{X}{|X|}$. In the proof of Theorem~\ref{thm2025-1-25-3}, we have $\n_{e_n}e_n = 0$, which means the integral curves of $e_n$
are geodesic. Also, by  Lemma~\ref{lem2025-1-23-2},  each hypersurface $\Sigma: |X| = c$ is totally umbilic in $(M, g)$. Thus, from a result in \cite{p-r}, $M$ is locally a twisted product $I\times_f \Sigma$ so that the metric tensor $g$ takes of the form
\be
g = dt^2 + f(t, x)^2 g_\Sigma,\label{eqn2025-5-18-4}
\ee
where $\frac{\partial}{\partial t} = \frac{X}{|X|} = e_n.$

Now let $(t, x^1, \cdots, x^{n-1})$ be a local coordinate system on $I \times \Sigma$. As in the proof of Theorem in \cite{chen}, $V = f\frac{\partial}{\partial t}$ is a torqued vector field with the torqued function $\vp = \frac{\partial f}{\partial t}$ and the torqued form $\o$ satisfies
$\o(\frac{\partial}{\partial t}) = 0$ and $\o(Y) = Y(\ln f)$ for $Y \perp \frac{\partial}{\partial t}$. In particular, the function $|V|^2 = f^2$ satisfies
$$
\n f^2 = 2f \n f
$$
and
$$
\n |V|^2 = 2\vp V + 2f^2 \o^\sharp.
$$
Thus, we obtain
\be
\n f= \frac{\vp}{f}V + f \o^\sharp = \vp \frac{\partial}{\partial t} + \sum_{i=1}^{n-1}\frac{\partial f}{\partial x^i}\frac{\partial}{\partial x^i}.\label{eqn2025-5-18-5}
\ee
On the other hand,  with respect to the coordinate system $(t, x^1, \cdots, x^{n-1})$ and the twisted product metric (\ref{eqn2025-5-18-4}), we have
\be
\n f = \frac{\partial f}{\partial t}\frac{\partial}{\partial t} + \frac{1}{f^2} \sum_{i=1}^{n-1}\frac{\partial f}{\partial x^i}\frac{\partial}{\partial x^i}.\label{eqn2025-5-18-6}
\ee
Comparing this to (\ref{eqn2025-5-18-5}), we have either $f \equiv 1$ and so $g$ is a (locally) product metric or $\frac{\partial f}{\partial x^i} = 0$ for all $i = 1, 2, \cdots, n-1$, which means that 
$f = f(t)$ and so the metric tensor $g$ is a (locally) warped product metric.

Finally, assume $f = 1$ on $M$ so that $M$ is a locally product $I \times \Sigma$. In this case, we hav e$X = |X|\frac{\partial}{\partial t}$ with respect to the coordinate $(t, x^1, \cdots, x^{n-1})$. Since $\n_{\partial _j} \partial_t = 0$, where $\partial_j = \frac{\partial}{\partial x^j}$ and
$\partial_t = \frac{\partial}{\partial t}$, we have
$$
\n_{\partial_j} X = \frac{\partial |X|}{\partial x^j} \frac{\partial}{\partial t}.
$$
On the other hand, since $\n_{\partial_j} X = \psi \frac{\partial}{\partial x^j}$, we obtain $\psi = R-\rho =0$, which means $X$is a Killing vector field.
\end{proof}


\bigskip
\noindent Seungsu Hwang\\
{Department of Mathematics }\\
Chung-Ang University \\
Heukseok-ro 84, Dongjak-gu, Seoul 06974, Korea \\
{\tt E-mail:seungsu@cau.ac.kr}\\

\bigskip
\noindent Gabjin Yun \\
Department of Mathematics and The Natural Science Research Institute\\
Myongji University \\
Myongji-ro 116, Cheoin-gu, Yongin, 17058, Korea \\
{\tt E-mail:gabjin@mju.ac.kr} \\

\end{document}